\documentclass[12pt]{amsart}
\usepackage[T1]{fontenc}
\usepackage[latin1]{inputenc}
\usepackage{amsmath}
\usepackage{amssymb}
\usepackage{graphics}
\newtheorem{theo}{\bf Theorem}[section]

\newtheorem{coro}[theo]{\bf Corollary}

\begin{document}

\title[The lamp lighting problem]{Note on the lamp lighting problem}
\author{Henrik Eriksson}
\address{NADA \\
   KTH \\
   SE-100 44 Stockholm, Sweden}
\email{henrik@nada.kth.se}
\author{Kimmo Eriksson}
\address{IMA, M{\"a}lardalens h{\"o}gskola \\ 
   Box 883 \\
   SE-721 23 V{\"a}ster{\aa}s, Sweden}
\email{kimmo.eriksson@mdh.se}
\author{Jonas Sj{\"o}strand}
\address{NADA \\ 
   KTH \\
   SE-100 44 Stockholm, Sweden}
\email{jonass@nada.kth.se}
\date{}

\begin{abstract}
We answer some questions concerning the so called $\sigma$-game of Sutner. 
It is played on a graph where each vertex has a lamp, the light of which 
is toggled by pressing any vertex with an edge directed to the lamp.

For example, we show that every configuration of lamps can be lit if and only 
if the number of complete matchings in the graph is odd. In the special case
of an orthogonal grid one gets a criterion for whether the number of 
monomer-dimer tilings of an $m\times n$ grid is odd or even.
\end{abstract}
\maketitle

\section{Introduction}
A janitor is walking through a museum in the morning, switching on the light 
in every room. The electric connections are a bit tricky, though. In every 
room is a button, but pressing this button toggles the light on/off not only
in the same room but also in all neighboring rooms! Can the janitor light up 
the entire museum?

An equivalent version of this problem was introduced by Sutner in 1989 
\cite{Su1}, where this game is called the $\sigma^+$-game on the graph of 
the museum. In the even more perverse $\sigma$-game, a button lights only 
its neighbors but not its own room. Such games have since been studied 
further by Sutner \cite{Su2,Su3}, Barua and Ramakrishnan \cite{BR} and
Goldwasser, Klostermeyer and Trapp \cite{GKT}. The commercially available 
"Lights out" game has the same rules and has been studied independently by 
Anderson and Feil \cite{AF} and, recently, by Dyrkacz, Eisenbud and  Maurer 
\cite{DEM}. Several precursors exist, such as the 
{\em Gale-Berlekamp switching game} from the sixties \cite{Sp}, where there 
are switches that toggle entire rows and columns, and the {\em Merlin game} 
studied by Pelletier \cite{Pe} and Stock \cite{St}, where the game board is 
no longer a grid. One should also note the similarity to the very general 
{\em lightbulb networks} of Kauffman \cite{Ka}, where the light condition of 
any lightbulb at a given time is determined by some rule (different for 
different lightbulbs) depending only on the light of the neighbors in the 
previous time-step.

In its most general setting, the $\sigma$-game is played on a directed graph 
where every vertex has a button and a lamp. At the beginning all lamps are 
switched off. When a button is pressed the lamps are toggled on/off on all 
vertices to which there is an edge from the pressed vertex. The objective is 
to light as many lamps as possible.

\section{When can all lamps be lit?}
Sutner \cite{Su1} studied, in particular, undirected graphs with a loop on every vertex (so that every button toggles its local lamp). He proved that in this special case it is always possible to light every lamp. The proof used linear algebra in a clever way, and Sutner remarked that no graph-theoretic proof was known. We will give such a proof, which also allows a generalization to a family of directed graphs.

We use the convention that a loop on a vertex $v$ contributes to both the out-degree and the in-degree of $v$.

\begin{theo}\label{th:jonasgeneral}
If $G$ is a directed graph on vertex set $V$, such that for each odd subset 
$U\subseteq V$ there is a vertex with odd out-degree in the induced subgraph 
on $U$, then it is possible to light all lamps.
\end{theo}
\begin{proof}
By induction over the number of vertices. The statement is trivially true for 
$|V|=1$. Suppose it holds for $|V|=n$ and consider the case $|V|=n+1$. 
If we choose any vertex $v$ and remove it and all its edges from $G$, then
the induction hypothesis still applies so it is possible to light all lamps. 
Call this an $n$-pressing with respect to $v$. If we 
apply this $n$-pressing to $G$ (that is, including $v$), then either $v$ too
lights up and we are finished, or $v$ is still dark. The only case left is 
when all $n+1$ such $n$-pressings leave one vertex dark.

Case 1: $n+1$ is even. Add all $n+1$ $n$-pressings. Every vertex will be 
lit an odd number of times and hence end up lit.

Case 2: $n+1$ is odd. Then there is a vertex $u$ with odd out-degree. Press 
the button at $u$. Let $U$ be the set of vertices now lit. Then $|U|$ is 
odd and $|V\backslash U|$ is even. Now add the $n$-pressings
with respect to each vertex in $V\backslash U$, one at a time. This will 
light the lamps in $V\backslash U$, while the lamps in $U$ will remain lit.
\end{proof}

We note that in order to satisfy the premise for single point subsets $U$, 
the graph must have a loop on each vertex . If all other edges come in 
opposite pairs (the undirected case), the number of edges in the subgraph 
induced by any odd subset $U$ is odd, so the premise will be satisfied in 
full and Sutner's result follows. 

The following corollary extends this result to graphs
with some of the edges bidirected (i.e.~opposite pairs) and some
uni-directed (i.e.~unpaired).

\begin{coro}\label{co:jonas}
If $G$ is a directed graph on vertex set $V$ with a loop on every vertex, such 
that the set of uni-directed edges form a complete bipartite graph on $V$, then
 it is possible to light all lamps.
\end{coro}
\begin{proof}
A complete bipartite graph on an odd vertex set $U$ has an even number of
edges, as one of the parts must be even. Thus the induced subgraph of $G$ 
on an odd subset $U\subseteq V$ has an odd number of edges (there is an
odd number of loops and the contribution from the bidirected edges is even)
and so some vertex must have odd out-degree. By Theorem \ref{th:jonasgeneral}, 
all lamps can be lit.
\end{proof}

If we take the complete bipartite graph in the corollary to be the edgeless graph $K_{0,|V|}$, then we obtain Sutner's result. For the next case, 
$K_{1,|V|\!-\!1}$, the museum interpretation might be the following. Some of
the exhibition rooms have been converted into offices. Buttons work as
before (own lamp and neighbors) but a new room, {\em the corridor}, is special.
The corridor button toggles all exhibition rooms (and the corridor) and
all office buttons toggle the corridor. Then all lamps can be lit! 

The parity argument used will of course work for any directed graph $G$ such 
that the induced subgraph of $G$ on any odd subset $U\subseteq V$ has an 
odd number of edges. The apparent greater generality of this statement
is illusory, for this property holds only for $G$ such that the corollary
applies (exercise!).

\section{How many lamps can always be lit?}
On a directed graph where there is at least one edge directed to each lamp, 
more than half the vertices can always be lit. This is a very easy exercise 
in the probabilistic method: Press each button with probability $0.5$. 
Then every lamp has probability one half of being lit, so the expected number 
of lit lamps is $|V|/2$.
Since at least one combination (when no button is pressed) has zero lit lamps, 
there must exist a combination with more than half the vertices lit. 

Spencer \cite{Sp} uses a not much more sophisticated argument to show that 
in the special case of the Gale-Berlekamp switching game, substantially more 
than half the lamps can be lit regardless of the initial configuration. No 
such thing holds for the general
$\sigma$-game as can be seen from the following construction. For a positive integer $n$, take $2^{n}-1$ lamps, labeled by the nonzero binary $n$-vectors. Now introduce $n$ buttons, such that button $i$ toggles the $2^{n-1}$ lamps that have a one in position $i$. Evidently every combination of pressing buttons (except for pressing none at all) lights $2^{n-1}$ lamps, so this is the maximal number.

\section{Directed games are equivalent to undirected games}
Although it would seem that the lamp-lighting game must be much more general on directed graphs than on undirected graphs, it turns out that in a very basic sense the games are equivalent: For each directed graph one can find an undirected graph on the same vertex set such that exactly the same configurations of lamps can be lit in the two games. The maximal number of lit lamps is related to the number of loops in the undirected graph.

\begin{theo}\label{th:jonas2}
Suppose that $k$ is the maximal number of lamps that can be lit in the game on a directed graph $G$.
Then there is an undirected graph $G'$ on the same vertex set, with $k$ loops, such that exactly the
same subsets of lamps can be lit.
\end{theo}
\begin{proof}
Let ${\bf A}=(a_{ij})$ be the adjacency matrix of $G$, that is, $a_{ij}=1$ if there is an edge from $i$ to $j$,
zero otherwise. Let ${\bf a}_i$ be the $i$th row of $\bf A$.
The starting position with all lights out is denoted by a vector of zeros, one for each lamp. Pressing vertex $i$ has the effect of adding (modulo 2) the corresponding row ${\bf a}_i$ to the position vector.
Hence the space of all lightable lamp configurations is the row space of $\bf A$ (modulo 2). 

An undirected graph $G'$ with $k$ loops has an adjacency matrix $\bf A'$ that is symmetric and has $k$ ones on the main diagonal. We must show that such a matrix $\bf A'$, with the same row space as $\bf A$, can be found. 

For convenience, renumber the lamps (the columns of $\bf A$) so that the first $k$ lamps constitute a maximal lightable subset of lamps. The row space of $\bf A$ is preserved under Gauss elimination to row echelon form. If necessary, renumber the $k$ first lamps so that, after elimination, the matrix has the block form
\begin{equation}\label{eq:mat}
\begin{pmatrix}
 \bf I & \bf B \\ 
 \bf 0 & \bf 0 \\
\end{pmatrix},
\end{equation}
where $\bf I$ is an identity matrix of size $r\le k$, and $\bf B$ is some arbitrary matrix. Now, we know that the vector ${\bf 1}_k=(1,1,\dots,0,0,\dots)$ with $k$ ones lies in the row space. The only way we can obtain the first $r$ ones in a linear combination of the rows is by adding all the $r$ nonzero rows, so this sum must be ${\bf 1}_k$. Hence each of the first $k-r$ columns of $\bf B$ has an odd number of ones, while the other columns have even numbers of ones. 

The matrix in (\ref{eq:mat}) is row equivalent to the symmetric matrix
$$\bf A' = \begin{pmatrix}
 \bf I & \bf B \\ 
 {\bf B}^{\rm T} & {\bf B}^{\rm T}\bf B \\
\end{pmatrix}.$$
A diagonal element of ${\bf B}^{\rm T}\bf B$ is the scalar product of the corresponding column of $\bf B$ with itself (mod 2).
Hence it is zero if the column has an even number of ones, and one otherwise. Therefore ${\bf B}^{\rm T}\bf B$ will have
$k-r$ ones on the diagonal, so $\bf A'$ is a symmetric matrix with $k$ ones on the diagonal.
\end{proof}

In particular, every directed graph on which all lamps can be lit is equivalent to some undirected
graph with a loop on every vertex.

\section{Pressing buttons in dark rooms only}
Now let us return to the unfortunate janitor. When touring the museum, he may feel uncomfortable pressing buttons in rooms where the light is already on. Is it possible for him to press buttons in dark rooms only? In the usual museum topology (an orthogonal grid) the answer is yes. More generally:

\begin{theo}\label{th:kimmo}
In each bipartite undirected graph with a loop on every vertex, one can light every lamp by pressing only vertices where the lamps are currently off.
\end{theo}
\begin{proof}
We know from Sutner's theorem that there exists a subset $V'$ of vertices such that if they are pressed, then all lamps will be lit. Let $X$ and $Y$ be the disjoint vertex sets in the bipartition of the graph. Begin by pressing the vertices in $V'\cap X$ in any order. Since there are no edges between vertices in $X$, the lamps on these vertices will all be off at the time they are pressed. We can now press the vertices in $V'\cap Y$ in any order. They must all be off, since we know that they will all be lit when we are done and they will only be toggled once in this process.
\end{proof}

A counter-example for arbitrary graphs looks as follows:
\begin{center}
\begin{picture}(40, 40)(0, -10)
\setlength{\unitlength}{0.6mm}
\put (0,0){\circle{3}}
\put (10,0){\circle{3}}
\put (20,0){\circle{3}}
\put (30,0){\circle{3}}
\put (15,10){\circle{3}}
\put (3,0){\line(1,0){4}}
\put (13,0){\line(1,0){4}}
\put (23,0){\line(1,0){4}}
\put (11,1){\line(1,2){3}}
\put (19,1){\line(-1,2){3}}
\end{picture}
\end{center}

The only way to light all lamps is to press all the vertices of the triangle, and none of the other two. However, as soon as we have pressed one vertex in the triangle all its three vertices will be lit, so the next we press must be a lit vertex.

\section{When can every lamp configuration be lit?}
A {\em complete matching} in a graph is a subset $M$ of the edges such that every vertex is incident to exactly one edge in $M$. 

\begin{theo}\label{th:jonas3}
Let $G$ be an undirected graph. Then every lamp configuration can be lit if and only if the number of complete matchings in $G$ is odd.
\end{theo}
\begin{proof}
Every lamp configuration can be lit if and only if the adjacency matrix $\bf A$ is invertible modulo 2. 
The determinant of $\bf A$ modulo 2 is the sum over all coverings of the vertices of $G$  by disjoint directed circuits. Since all circuits of length $\ge 3$ come in pairs (two directions), they give zero contribution modulo 2. It remains coverings of all vertices by disjoint circuits of length 2 (undirected edges between pairs of vertices) and length 1 (loops), that is, precisely complete matchings. The adjacency matrix $\bf A$ is invertible if and only if its determinant is nonzero modulo 2, which is the case if the number of complete matchings is odd. 
\end{proof}

Barua and Ramakrishnan \cite{BR} and Goldwasser, Klostermeyer and Trapp \cite{GKT} have determined, by  methods entirely different from ours, for which sizes $m \times n$ of orthogonal grids (with loops on every vertex) that every lamp configuration can be lit. Their answer is: if and only if $p_m(\lambda)$ and $p_n(1+\lambda)$ are relatively prime,
where $p_m(\lambda)$ is the binary Chebyshev polynomial defined by the recurrence:
\begin{equation}\label{eq:BR}
  p_m(\lambda)=\lambda p_{m-1}(\lambda)+p_{m-2}(\lambda) \mbox{ for } m\ge2\mbox{,\ \ } 
  p_1(\lambda)=\lambda\mbox{,\ \ }p_0(\lambda)=1.
\end{equation}
By Theorem \ref{th:jonas3}, the result of Barua and Ramakrishnan is a 
statement about the parity of the number of complete matchings in a 
square grid with a loop on every vertex. But a complete matching in 
such a graph is precisely what is called a monomer-dimer tiling of the 
square grid; the loops are monomers and  edges covering two vertices are dimers. The number of monomer-dimer tilings of the $m\times n$ grid is a famous open problem, cf.\@ Finch's webpage \cite{SF}. Barua and Ramakrishnan's result and our Theorem \ref{th:jonas3} combine to the following partial result on the monomer-dimer problem.

\begin{coro} The number of monomer-dimer tilings of the $m\times n$  grid is odd if and only if $p_m(\lambda)$ and $p_n(1+\lambda)$ are relatively prime modulo 2.
\end{coro}

\section{How few lamps can be lit on an infinite grid}
We conclude by finding the minimum number of lit lamps on the infinite 
orthogonal grid with a loop on each vertex.

Pressing one vertex (at the origin, say) lights five lamps. 
By pressing some neighbour vertices one can of course switch them off again,
but in this process new lamps are lit. Is it possible that
pressing a complicated vertex pattern might result in fewer than five
lit lamps? 

We will prove that five is in fact the minimum and that the
only way to light exactly five lamps is by pressing a {\em mikado diamond} (see, and listen to, \cite{mikado}).
These are diamond-shaped dot patterns that leave only the center lamp and
the four extreme lamps lit; the first mikado diamond is the single dot
and each subsequent diamond doubles the size of the previous one. 
\begin{figure}[hb]
  \begin{center}
{
\setlength{\unitlength}{0.25mm}
\begin{picture}(60, 150)(-20, -70)
\put (0,0){\circle*{5}}
\end{picture}
\begin{picture}(60, 150)(5, -80)
\put (0,0){\circle*{5}}
\put (-10,-10){\circle*{5}}
\put (0,-10){\circle*{5}}
\put (10,-10){\circle*{5}}
\put (0,-20){\circle*{5}}
\end{picture}
\begin{picture}(130,150)(0, -100)
\put (0,0){\circle*{5}}
\put (-10,-10){\circle*{5}}
\put (0,-10){\circle*{5}}
\put (10,-10){\circle*{5}}
\put (-20,-20){\circle*{5}}
\put (20,-20){\circle*{5}}
\put (-30,-30){\circle*{5}}
\put (-20,-30){\circle*{5}}
\put (0,-30){\circle*{5}}
\put (20,-30){\circle*{5}}
\put (30,-30){\circle*{5}}
\put (0,-60){\circle*{5}}
\put (-10,-50){\circle*{5}}
\put (0,-50){\circle*{5}}
\put (10,-50){\circle*{5}}
\put (-20,-40){\circle*{5}}
\put (20,-40){\circle*{5}}
\end{picture}
\begin{picture}(160, 150)(0, -140)
\put (0,0){\circle*{5}}
\put (-10,-10){\circle*{5}}
\put (0,-10){\circle*{5}}
\put (10,-10){\circle*{5}}
\put (-20,-20){\circle*{5}}
\put (20,-20){\circle*{5}}
\put (-30,-30){\circle*{5}}
\put (-20,-30){\circle*{5}}
\put (0,-30){\circle*{5}}
\put (20,-30){\circle*{5}}
\put (30,-30){\circle*{5}}
\put (-40,-40){\circle*{5}}
\put (-20,-40){\circle*{5}}
\put (0,-40){\circle*{5}}
\put (20,-40){\circle*{5}}
\put (40,-40){\circle*{5}}
\put (-50,-50){\circle*{5}}
\put (-40,-50){\circle*{5}}
\put (-30,-50){\circle*{5}}
\put (30,-50){\circle*{5}}
\put (40,-50){\circle*{5}}
\put (50,-50){\circle*{5}}
\put (-60,-60){\circle*{5}}
\put (0,-60){\circle*{5}}
\put (60,-60){\circle*{5}}
\put (-70,-70){\circle*{5}}
\put (-60,-70){\circle*{5}}
\put (-40,-70){\circle*{5}}
\put (-30,-70){\circle*{5}}
\put (-10,-70){\circle*{5}}
\put (0,-70){\circle*{5}}
\put (10,-70){\circle*{5}}
\put (30,-70){\circle*{5}}
\put (40,-70){\circle*{5}}
\put (60,-70){\circle*{5}}
\put (70,-70){\circle*{5}}
\put (0,-140){\circle*{5}}
\put (-10,-130){\circle*{5}}
\put (0,-130){\circle*{5}}
\put (10,-130){\circle*{5}}
\put (-20,-120){\circle*{5}}
\put (20,-120){\circle*{5}}
\put (-30,-110){\circle*{5}}
\put (-20,-110){\circle*{5}}
\put (0,-110){\circle*{5}}
\put (20,-110){\circle*{5}}
\put (30,-110){\circle*{5}}
\put (-40,-100){\circle*{5}}
\put (-20,-100){\circle*{5}}
\put (0,-100){\circle*{5}}
\put (20,-100){\circle*{5}}
\put (40,-100){\circle*{5}}
\put (-50,-90){\circle*{5}}
\put (-40,-90){\circle*{5}}
\put (-30,-90){\circle*{5}}
\put (30,-90){\circle*{5}}
\put (40,-90){\circle*{5}}
\put (50,-90){\circle*{5}}
\put (-60,-80){\circle*{5}}
\put (0,-80){\circle*{5}}
\put (60,-80){\circle*{5}}
\end{picture}
}
{\setlength{\unitlength}{1mm}
\begin{picture}(10, 30)(0, -17.5)
\put (-20,-1){$\implies$}
\put (0,0){\circle{2}}
\put (10,0){\circle{2}}
\put (-10,0){\circle{2}}
\put (0,10){\circle{2}}
\put (0,-10){\circle{2}}
\put (0,-9){\line(1,1){9}}
\put (-9,0){\line(1,1){9}}
\put (-9,0){\line(1,-1){9}}
\put (0,9){\line(1,-1){9}}

\end{picture}
}
    \caption{The mikado diamonds light five lamps}
    \label{fig:mikados}
  \end{center}
\end{figure}
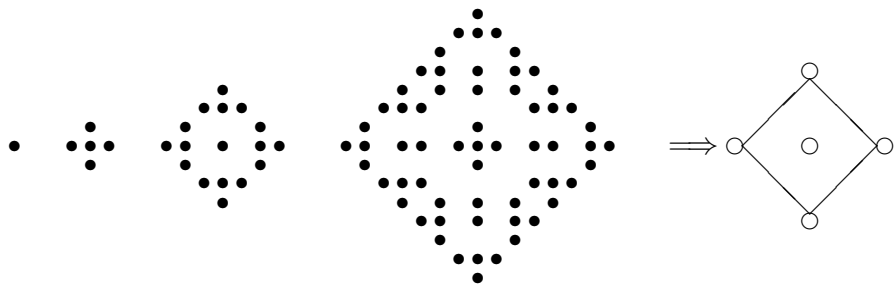

The five lamps property follows from the recursive construction indicated in
the next figure. Each mikado pattern is obtained by superposition of five 
mikado patterns of the next smaller size. It is clear that in this way
we get another pattern leaving only the center lamp and the extreme lamps lit,
for each the other eight lamps involved is lit twice, i.e. not at all.

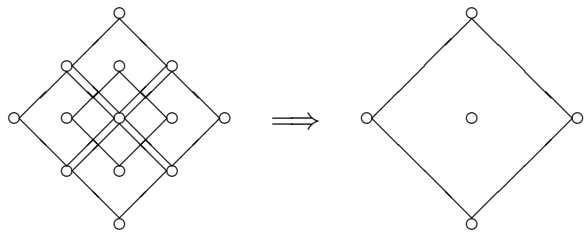
\begin{figure}[hb]
  \begin{center}
{\setlength{\unitlength}{0.7mm}
\begin{picture}(65, 40)(-20, -20)
\put (0,0){\circle{2}}
\put (10,0){\circle{2}}
\put (20,0){\circle{2}}
\put (-10,0){\circle{2}}
\put (-20,0){\circle{2}}
\put (10,-10){\circle{2}}
\put (10,10){\circle{2}}
\put (-10,-10){\circle{2}}
\put (-10,10){\circle{2}}
\put (0,10){\circle{2}}
\put (0,-10){\circle{2}}
\put (0,20){\circle{2}}
\put (0,-20){\circle{2}}
\put (0,-19){\line(1,1){9}}
\put (0,1){\line(1,1){9}}
\put (0,-9){\line(1,1){9}}
\put (-10,-9){\line(1,1){9}}
\put (10,-9){\line(1,1){9}}
\put (1,0){\line(1,1){9}}
\put (-9,0){\line(1,1){9}}
\put (-19,0){\line(1,1){9}}
\put (-9,10){\line(1,1){9}}
\put (-9,-10){\line(1,1){9}}
\put (1,0){\line(1,-1){9}}
\put (-9,0){\line(1,-1){9}}
\put (-19,0){\line(1,-1){9}}
\put (-9,10){\line(1,-1){9}}
\put (-9,-10){\line(1,-1){9}}
\put (0,-1){\line(1,-1){9}}
\put (-10,9){\line(1,-1){9}}
\put (0,9){\line(1,-1){9}}
\put (0,19){\line(1,-1){9}}
\put (10,9){\line(1,-1){9}}

\end{picture}
\begin{picture}(40, 40)(-20, -20)
\put (-40,-2){$\implies$}
\put (0,0){\circle{2}}
\put (20,0){\circle{2}}
\put (-20,0){\circle{2}}
\put (0,20){\circle{2}}
\put (0,-20){\circle{2}}
\put (0,-19){\line(1,1){19}}
\put (-19,0){\line(1,1){19}}
\put (-19,0){\line(1,-1){19}}
\put (0,19){\line(1,-1){19}}
\end{picture}
}
\caption{Superposition of five mikado diamonds makes a larger one}
\label{fig:superposition}
\end{center}
\end{figure}

Unexpectedly, one retrieves the smaller diamond by 
erasing every other row and every other column in the larger one. 
This is a consequence of the linear relation stating that the 
lamp at $(i,j)$ is dark:
$$x_{i,j}+ x_{i,j+1}+ x_{i,j-1}+ x_{i+1,j}+ x_{i-1,j}
\equiv 0 \pmod{ 2},$$ 
where $x$ is the button pressing matrix. Where darkness rules, 
many such relations hold, and by adding the five relations
belonging to lamp $(i,j)$ and its neighbors we obtain
$$x_{i,j}+ x_{i,j+2}+ x_{i,j-2}+ x_{i+2,j}+ x_{i-2,j}
\equiv 0 \pmod{ 2},$$ 
proving that darkness will prevail after the erasing process.
By our next theorem, we must have retrieved the mikado diamond.

The distance between the leftmost and the rightmost lamp is $2^k$ for
the $k$-th diamond. It is clear from the stated erasing property
that increasing $k$ just means adding more detail to the same picture.
For large $k$, the pattern is a fractal
with a recurrent mikado-like figure appearing in all sizes and 
orientations.   

\begin{figure}[hb]
  \begin{center}
\resizebox{110mm}{!}{\includegraphics{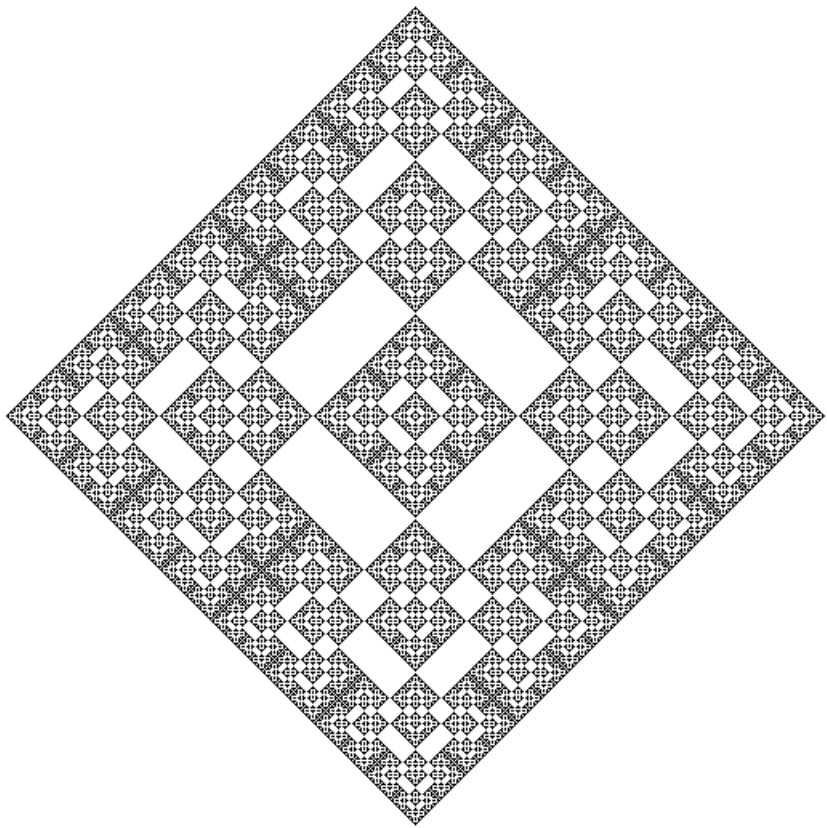}}
    \caption{The mikado pattern. (Can you see him?)}
    \label{fig:Mikado}
  \end{center}
\end{figure}

\begin{theo}
  No button pressing pattern will light one, two, three or four lamps
  and the only patterns that light five lamps are the mikado diamonds. 
\end{theo}
\begin{proof}
  
To see that a set of pressed vertices can never light one, two or three 
lamps, consider the smallest rectangle (with horizontal and vertical sides)
 covering such a set. It is clear that at least one lamp just outside each of the 
four sides of the rectangle must be lit. Any pattern that lights these four
lamps only must have $45\deg$ sloping boundaries, much the same as a mikado
pattern but it is not immediately clear why all four boundaries should have 
size $2^k$. However, considering the alternating pattern in the subdiagonal 
just inside the diagonal boundary, we can conclude that the lamp distances
are at least even numbers. 
If we apply the erasing process to the pattern, we obtain a smaller pattern
with exactly the same property, contradicting the assumed minimality. 

The same reasoning applies to the five lamps case also. First, note that the 
whole pattern must be symmetric, otherwise superposing its mirror image would
cancel the four boundary lamps and leave at most two interior lamps lit. So the
fifth lit lamp must be at the center and we can perform the erasing procedure
until we are down to the smallest one-button mikado pattern. 
\end{proof}

A final observation is that except for 1, 2, 3 and 4, any desired number $r$ of
lamps can be lit by pressing $r-4$ buttons that are diagonally consecutive.

\end{document}